\documentstyle [12pt]{article}

\begin{document}
\begin{center}
\LARGE {\bf The Length of the Shortest Closed Geodesics\\
       on a Positively Curved Manifold}\vspace{5mm}\\
\large Yoe Itokawa \normalsize and \large Ryoichi Kobayashi
\end{center}
\def\thefootnote{\fnsymbol{footnote}}
\footnotetext[0]{2000{\it Mathematics Subject Classification.}
 53C20, 53C22.}
\footnotetext[0]{{\it Key words and phrases.}
 Riemannian manifold, positive sectional curvature, closed 
 geodesics, loop spaces, Morse theory.}
\begin{abstract}
 We give a metric characterization of the Euclidean sphere in 
 terms of the lower bound of the sectional curvature and the 
 length of the shortest closed geodesics.
\end{abstract}

\normalsize
\newcommand{\eu}{\mbox{\bf R}}
\newcommand{\z}{\mbox{\bf Z}}
\newcommand{\vnn}{\vspace{5mm} \par \noindent}
\newcommand{\svnn}{\vspace{1mm}\newline}
\newcommand{\nln}{\newline\noindent}
\newcommand{\n}{\mid}
\newcommand{\wt}{\widetilde}
\newcommand{\ep}{\varepsilon}

\section{Introduction}

Let $M$ be a complete connected Riemannian manifold of dimension 
$d$ and class $C^{\infty}$.  The study of global structure of 
closed geodesics on $M$ {\it vis a vis} certain quantitative 
restrictions on the sectional curvature $K$ of $M$ has attracted 
considerable interest among researchers.  Henceforth, we assume 
$k$ to be a positive constant.  It follows straightforwardly 
from Morse-Schoenberg index comparison that if $K \geq  k^2$ on 
all tangent 2-planes of $M$, then there must exist on $M$ a closed 
geodesic whose length is $\leq 2\pi/k$ representing the lowest 
dimensional homology class of $M$.  The purpose of the present 
paper is to describe a rigidity phenomenon observed when this 
length is extremal on $M$.  More precisely, we prove 

\vnn
{\bf Main Theorem.}\quad{\it
 If $M$ satisfies $K \geq k^2$ and 
 if the shortest closed geodesics on $M$ have length $=2\pi/k$, 
 then $M$ is isometric to $S^d_k$, the Euclidean sphere of radius 
 $1/k$ in $\eu^{d+1}$.
}
\vspace{5mm}

Note that we make no assumption about the geodesics' having no 
self-intersections.  There exists an example of a 2-dimensional 
smooth surface all of whose shortest closed geodesics have 
self-intersections.  These examples have some regions where the 
curvature is negative.  Calabi and Cao [CC] has proved that on a 
positively curved surface, at least one of the shortest closed 
geodesics is always without self-intersections.
\footnotetext{{\it Acknowledgement.}
 The authors wish to acknowledge their 
 gratitude to D. Gromoll, G. Thorbergsson and W. Ballman for 
 providing them with valuable suggestions and information. 
}

We now mention some related rigidity phenomena.  Previously, 
Sugimoto [Su], improving on an earlier work of Tsukamoto [Ts], 
proved

\vnn
{\bf Theorem A.}\quad{\it
 Suppose that $M$ satisfies $4k^2 \geq K \geq k^2$.  If $d$ is 
 odd, assume that $M$ is simply connected.  Then, if $M$ has a 
 closed geodesic of length $2\pi/k$, it is isometric to $S^d_k$.
}
\vspace{5mm}

Recall that under the curvature assumption of Theorem A, if $M$ 
is simply connected, the celebrated Injectivity radius theorem, 
which is primarily due to Klingenberg (see [CE (\S\S 5.9,10)], 
[GKM, \S\S7.5,7] and also [CG], [K] and [Sa2]) states that all 
closed geodesics on $M$ have length $\geq \pi/k$.

However, we point out that, in general, an assumption on the 
length of the shortest closed geodesic is a nontrivially weaker 
condition than an upper bound on the sectional curvature.  In 
fact, it is possible to construct, for any given $k$ and $\delta$, 
a Riemannian metric on $S^2$ with $K \geq k^2$ and the length of 
the shortest closed geodesic $\delta$-close to $2\pi/k$ but whose 
curvature grows arbitrarily large ($S^2$ like surface with highly curved 
``equator").  This construction means that, 
from the viewpoint of rigidity theorems in Riemannian geometry, 
imposing an upper bound on the curvature is not natural in 
characterizing a Euclidean sphere among complete Riemannian 
manifolds with $K \geq k^2$ having a shortest closed geodesic of 
length just $2\pi/k$.  For more informations on 
the curvature bounds and the lengths of closed geodesics, see, 
for instance, [Sa2].  

In the spacial case of dimension 2, we have

\vnn
{\bf Theorem B} ({\it Toponogov} [T]).\quad{\it
 Suppose that $M$ is an abstract surface with Gauss curvature 
 $K \geq k^2$.  If there exists on $M$ a closed geodesic without 
 self-intersections whose length $=2\pi/k$, then $M$ is isometric 
 to $S^2_k$.
}
\vspace{5mm}

However, in higher dimensions, there are lens spaces of constant 
sectional curvature $k^2$ whose geodesics are all closed with 
the prime ones have no self-intersections, and they are either

(a) homotopic to 0 and have length  $=2\pi/k$, or

(b) homotopically nontrivial and can be arbitrarily short.

\noindent
{\it See }[Sa1]
\vspace{5mm}

Of course, it follows from our Main Theorem that

\vnn
{\bf Corollary.}\quad{\it
 If $K \geq k^2$ and the shortest closed geodesics that are 
 homotopic to 0 in $M$ have the length $2\pi/k$, then the universal 
 covering of $M$ must be isometric to $S^d_k$.
}

Note also that Theorem B is false without the assumption that 
the closed geodesics have no self-intersections.  In fact, for 
any $k$, one can construct an ellipsoid in $\eu^3$ which possesses 
a prime closed geodesic of length $= 2\pi/k$ and whose curvature 
is $> k^2$.

Finally, we mention a previous related result of the first author 
which gives another rigidity solution for the nonsimply connected 
case.

\vnn
{\bf Theorem C} ({\it Itokawa} [I1,2]).\quad{\it
 If the Ricci curvature of $M$ is $\geq (d-1)k^2$ and if the 
 shortest closed geodesics on $M$ have the length $\geq \pi/k$, 
 then either $M$ is simply connected or else $M$ is isometric 
 with the real projective space all of whose prime closed geodesics 
 have length $= \pi/k$.
}
\vspace{5mm}

It is not yet known if our Main Theorem remains true if the 
assumption on the sectional curvature is weakened to that on the 
Ricci curvature.  However, we point out that examples were shown 
in [I1,2] so that for the Ricci curvature assumption, the shortest 
closed geodesics may have length arbitrarily close to $2\pi/k$ 
without the manifold's even being homeomorphic to $S^d$.  This 
indicates how delicate the Ricci curvature assumption could be.
\vspace{5mm}

 \section{Preliminaries}

The purpose of this section is to collect together all the 
well-known facts and results which will be used in proving the 
Main Theorem and also to set straight our notational conventions 
and normalizations.  In this paper, we agree that by the term 
{\it curve} we mean an absolutely continuous mapping $c : \eu
 \longrightarrow M$ whose derivative is defined alomost everywhere 
and is $L_2$ on each closed interval.  We refer to the restriction 
of a curve to any closed interval as an {\it arc}.  If $c$ is a 
curve and $a < b$ are reals, we write $c_{a,b}$ to denote the arc
 $c \n_{[a,b]}$.  If $c$ happens to be differentiable, the normal 
bundle; respectively, the unit normal bundle of $c$; which are 
in fact bundles over $\eu$, are denoted $\perp c$; respectively, 
$U \perp c$.  We shall call a curve $c$ closed if $c(s+1) = c(s)$ 
for all $s$.  We denote the set of all absolutely continuous closed 
curves with $L^2$ derivative in $M$ by $\Omega$.
\vspace{5mm}

For fixed $a$, $b$, let ${\cal A}_{a,b}$ denote the set of all 
arcs $[a,b] \longrightarrow M$.  It is known that ${\cal A}_{a,b}$ 
has the structure of a Riemannian Hilbert manifold where the inner 
product is given by the natural $L_2$ inner product of variation 
vector fields along a curve.  The restriction map $c \longmapsto 
 c_{0,1}$ embeds $\Omega$ in ${\cal A}_{0,1}$ as a closed submanifold 
so that the Riemannian structure pulls back on $\Omega$.  We refer 
to the paper [GM2], the book [K], and references sited therein for 
details.
\vspace{5mm}

For $\gamma \in {\cal A}_{a,b}$, we define the space $
 {\cal V}'_{\gamma}$ of all square integrable vector fields $v \in 
 T_{\gamma}{\cal A}_{a,b}$ along $\gamma$ such that $v(a) = 0$ and 
$v(b) = 0$ and such that $v(s) \in \perp_s \gamma$ for all $s$, 
whereever $\gamma$ is differentiable.  If $c \in \Omega$, we also 
define the space ${\cal V}_c$ of all $v \in T_c\Omega$ with $v(s)
 \in \perp_s c$ almost everywhere.  Then, ${\cal V}'_{c_{0,1}}$ 
is canonically embedded in ${\cal V}_c$.
\vspace{5mm}

We normalize the energy of $\gamma \in {\cal A}_{a,b}$ by
$$
 E(\gamma) := \int_a^b \n \gamma'(s)\n^2ds \,\, .
$$
Also, we denote by $L(\gamma)$ the length of $\gamma$ in the usual 
sense.  
Thus, in our convention, $L(\gamma)^2 \leq (b - a) E(\gamma)$ 
with equality if and only if $\gamma$ is parametrized proportinal 
to arclength.  The term {\it geodesic} is always understood to 
mean a nonconstant geodesic.  For $u \in UTM$, the unit tangent 
bundle, we denote by $c_u$ the geodesic $s \mapsto \exp su$.  
Recall that the critical points of $E$ on $\Omega$ are closed 
geodesics and the constant curves.  See either [CE], [GKM], or 
[Mi].
\vspace{5mm}

Let $c$ be a geodesic and $a < b \in \eu$.  The Hessian of $E$ 
at $c_{a,b}$, here regarded as a symmetric bilinear form on 
$T{\cal A}_{a,b}$, is denoted $H_a^b$.  We remind the reader that 
if $v \in {\cal V}'_{c_{a,b}}$ and is differentiable outside of 
finitely many points, or if $c \in \Omega$, $a = 0$, $b = 1$, 
and $v \in {\cal V}_c$ is differentiable outside of finitely many 
points in $(a,b)$, then $H^b_a(v,v)$ is given by the {\it index 
integral}
$$
 - 2 \left \{ \int_a^b (\langle v''(s) , v(s) \rangle +
  \n v(s)\n^2 \n c'(s)\n^2 K_{v(s) \wedge c'(s)} ds +
  \sum_s \langle v(s) , \Delta_s v'(s) \rangle \right \}
$$
where 
$$
 \Delta_s v'(s) = v'(s_+) - v'(s_-)
$$
denotes the jump in $v'(s)$ at one of its finitely many points 
of discontinuity in the open interval $(a,b)$.  See, for example, 
[Mi], [Bo], or [BTZ]. We write $\iota'(c_{a,b})$ to denote the 
index of $H_a^b\n{\cal V}'_{c_{a,b}}$.  If $c$ is closed, we put 
$H := H_0^1\n{\cal V}_c$ and $\iota(c)$ its index.  We recall the
basic inequality
$$
 \iota(c) \geq \iota'(c) = \sum_{0 < s < 1} \nu'(c_{0,s})
$$
where $\nu'(c_{0,s})$ is the dimension of the space of Jacobi 
fields in ${\cal V}'_{c_{0,s}}$.  In this notation, we state the 
following well-known theorem, which is primarily due to Fet [F].

\vnn
{\bf Theorem D.}\quad{\it 
 Assume that $M$ satisfies $K \geq k^2$.  Then there exists a 
 closed geodesic $c$ on $M$ such that $L(c) \leq 2\pi/k$ and 
 $\iota(c) \leq d-1$.
}
\vspace{5mm}

For each $r \in \eu$, we denote by $\Omega^r$ ({\it respectively}, 
$\Omega^{=r}$ and $\Omega^{<r}$) the subspaces $\{ c \in \Omega
 \, : \, E(c) \leq r \quad \mbox{({\it respectively}, $= r$ and 
 $< r$)} \}$.  However, $\Omega^0=\Omega^{=0}$ is identified with 
$M$ itself and so denoted also by $M$.  It is well-known that 
the energy functional $E$ sataisfies the famous Condition C of 
Palais and Smale.  See, for example [GM2].  The significance of 
this for us is that, as far as global variational-theoretic 
properites are concerned, we can treat $E$ as if it were a proper 
function defined on a locally compact manifold.

Alternatively, we can work on the finite-dimensional approximation 
of $\Omega^r$ {\it \`{a} la} [Mi (\S 16)] or [Bo].  While this 
has the advantage that it simplifies the analytical aspect of the 
argument, we prefer to use the infinite-dimensional argument in 
\S3 because of the ease by which we can write the variational vector 
fields explicitly.  Not that we could write the corresponding 
fileds explicitly in the finite-dimensional approximation, but 
the actual expressions would be unpleasantly complicated. 

We must later consider a more general functional $F$ on $\Omega$ or 
$'\Omega$ (in \S4). 
Let $c \in \Omega$ be a critical point of $F$.  Then $T_c \Omega$ 
decomposes into a direct sum
$$
 T_c \Omega = {\cal P} \oplus {\cal N} \oplus {\cal Z}
$$
where $\cal P$, $\cal N$ and $\cal Z$ are the spaces on which 
the Hessian $H_F$ of $F$ at $c$ is positive definite, negative 
definite and zero respectively.  We write $\| \cdot \|$ for the 
norm in $T_c \Omega$.  Then, we can state the following important 
fact due to Gromoll and Meyer [GM1] (See especially the note on
{\it p}.362).

\vnn
{\bf Theorem E} (Generalized Morse Lemma).\quad{\it
 In the setting described above, there exists a neighborhood $U$ 
 of $c$, a coordinate chart
 $$
  \xi_c \, : \, U \longrightarrow T_c \Omega \, ,
 $$
 with respect to which $F$ takes the form
 $$
  F \circ \xi_c^{-1}(v) = \| x \|^2 - \| y \|^2 + f(z) + F(c)
 $$
 where $x$, $y$ and $z$ are the orthogonal projections of $v \in
  \xi_c^{-1}(U)$ on $\cal P$, $\cal N$ and $\cal Z$ respectively, 
 and $f$ is a function whose Taylor series expansion at $z = 0$ 
 starts with the term of degree at least 3 in $z$ or  equivalently 
 with vanishing Hessian.  For this  decomposition, $c$ need not 
 be an isolated critical point of $F$, but if $F$ has other 
 critical points in $U$, their images in $T_c \Omega$ are all 
 contained in $\cal Z$.
}
\vspace{5mm}

The chart $\xi_c$ is often called the {\it Gromoll-Meyer-Morse 
chart} at $c$ with respect to $F$.

We put $U_c^- := \xi_c^{-1}({\cal N})$ and $U_c^{-0} := 
 \xi_c^{-1}({\cal N} \oplus {\cal Z})$ and call them the {\it 
strong unstable submanifold} and the {\it weak unstable submanifold} 
of $F$ at $c$ respectivley, even though we make no assumption 
that $\dim \cal Z$ is finite.

Suppose that $a \in \eu$.  We set $\Omega^a_F := \{c \in \Omega
 \, : \, F(c) \le a \}$.  Let $I$ be the interval $[-1,1]$. 
Suppose $c$ is a critical point of $F$ with $a := F(c)$ and 
$\iota := \mbox{index} \,\, H_{F \n c} = \dim \cal N$.  Let $U$ 
be a neighborhood of $c$ as defined in Theorem E. 

We call a differentiable embedding 
$\sigma:(I^{\iota},\partial I^{\iota}) \longrightarrow 
(\Omega,\Omega^a_F-U)$ a {\it weak unstable simplex} 
(resp. {\it strong unstable simplex}), if there exists 
a smaller neighborhood $W$ of $c$, $c\in W\subset U$, 
so that $\sigma(I^{\iota}) \cap W$ coincides with $\xi^{-1}(\cal N\oplus \cal Z)$ 
(resp.$\xi^{-1}({\cal N}) \cap W$). 
If $\sigma$ is a weak unstable simplex of $F$ at $c$ in this sense, it is clear that 
$\sigma \cap W$ must be contained in the topological cone
$$
 \{ \gamma \in W \, ; \, 
  H_F(\xi(\gamma),\xi(\gamma)) \leq 0 \}
$$
containing $\xi^{-1}(N)$. Distinguishing from weak or strong unstable simplex 
just defined above, 
we mean, by {\it unstable simplex} of $F$ at $c$, a differentiable embedding 
$\sigma:(I^{\iota},\partial I^{\iota}) \longrightarrow (\Omega,\Omega^a_F-U)$ 
such that $\sigma(0)=c$ and $F|_{\sigma}\leq a$. 

We shall say that a critical point $c$ of $F$ is {\it nondegenerate} 
if ${\cal Z} = \{ 0 \}$.  Note that this agreement is different 
from the often-used convention of calling a closed geodesic 
nondegenerate if $\cal Z$ is the $S^1$ orbit of the geodesic. 
With our convention, a closed geodesic is never a nondegenerate 
critical point for $E$ because of the $S^1-$action.  
We put $a:= F(c)$ and write 
$\Omega^r := \{ \gamma \in \Omega \, : \, 
 F(\gamma) \leq r\}$.  If $c$ is a nondegenerate critical point 
of $F$, then of course $c$ is an isolated critical point and, 
for some $\ep > 0$, the strong unstable simplexes at $c$ 
represent a nontrivial class in the relative homotopy group 
$\pi_{\iota}(X^{a+\ep},X^{a-\ep})$. 
\vspace{5mm}

\section{Proof of the Main Theorem}

It is
clear that, in order to prove Theorem 1, we need to consider only one
specific $k$.  So, hereafter we assume that $M$ satisfies $K \geq
k^2$ where $k:=2\pi$.  In the present section, we further assume
that $M$ contains no closed geodesic of length $<1$, or
equivalently that there are no critical points of $E$ in
$\Omega^{<1}-M$.  It now remains for us to prove that then $M$ is
isometric to $S^d_{2\pi}$.\vnn \indent We set\vspace{5mm}
\begin{center}$\cal C$ $:=\{c \in \Omega\,;\,c$ is a closed
geodesic of length 1 and $\iota(c)=d-1\}$ \end{center}\vspace{5mm}
and\vspace{5mm} \begin{center} ${\cal C}^*$ $:=\{c \in \cal C\,;\,\,$
an unstable simplex of $E$ at $c$ represents\\ a nontrivial element
in $\pi_{d-1}(\Omega,M)\}\,\,.$ \end{center}\vspace{5mm}
Theorem D and the Morse-Schoenberg index comparison 
assert that $\cal C \not=\emptyset$.  
If we can asume that each $c \in \cal C$ has an isolated 
critical $S^1$-orbit, the technique of Gromoll and Meyer [GM2] 
fairly readily shows that ${\cal C}^*$ too is non-empty.  
In our case, however, it will be precisely 
one of our points that no $c \in {\cal C}^*$ has an isolated critical orbit. 
Under the stronger hypothesis of $4k^2 \geq K \geq k^2$, 
Ballman [Ba] showed that all closed geodesics have nontrivial 
unstable simplexes.  However, he makes essential use of the 
upper bound for $K$ which is not available to us. 
Nonetheless, we shall still prove in the next section,

\vnn {\bf Lemma 1.}\quad{\it
 Under the assumptions of this section, ${\cal C}^*$ is nonempty
 and is a closed set in $\Omega$.
}
\vspace{5mm}

In this section, we accept Lemma 1 for the time being, and prove

\vnn {\bf Lemma 2.}\quad{\it
 For each 
 $c \in {\cal C}^*$, there is a neighborhood $\cal U$ of $c'(0)$ in 
 $UT_{c(0)}M$ such that whenever $u \in \cal U$ and $\tau$ is any 
 tangent 2-plane containing $c_u'(s)$ for some $s \in \eu$, then 
 $K(\tau)=k^2$ ($c_u$ being 
 the geodesic determined by the initial condition $c_u'(0)=u$).
}
\vspace{5mm}

We prove Lemma 2 by proving a sequence of other Lemmas (from 3 to 8). 
The idea for proving Lemma 2 is to construct for every $c \in \cal C$ a 
specific unstable simplex $\tau$ which is homotopic to the strong 
unstable simplex and a deformation of such a $\tau$ so that, unless 
the conclusion of the lemma is met, $\tau$ is deformed into 
$M \subset \Omega$, which is a contradiction if $c \in {\cal C}^*$.  
First, we show

\vnn
{\bf Lemma 3.}\quad{\it 
 If $c \in \cal C$,
 then for any $s \in \eu$ and any $v \in \perp_sc$, $
  K(c'(s) \wedge v) \equiv k^2$.
}
\vspace{5mm}

{\it Proof.} 
 Assume that, for some $s_1 \in \eu$ and $v_1 \in \perp_{s_1}c$,
 $K(c'(s_1) \wedge v_1) > k^2$.  By virtue of the natural $S^1$-action on
$\Omega$, it is no loss of generality to assume that $0 < s_1 < 1/2$.  Now, we
define a real number $\delta$ as follows.  If there is a point in $(0,s_1]$
which is conjugate to 0 along $c$, we choose any $\delta$ so that
$s_1<\frac12-\delta<\frac12$.  If, on the other hand, there is no conjugate
point in $(0,s_1]$, there is a unique Jacobi field $Y$ along $c$ with $Y(0)=0$
and $Y(s_1)=v_1$, and by a consequence of the original Rauch comparison 
theorem [CE (\S 1.10, Remark, p.35)], there is an $s_2$, $s_1<s_2<1/2$ 
so that $Y(s_2)=0$.  In this case, we choose $\delta$ so that 
$s_2 < \frac12 - \delta < \frac12$.  In either case, we have 
$\iota'(c_{0,\frac12-\delta}) \geq 1$.  On the other hand, by the 
Morse-Schoenberg index comparison with $S^d_k$, we have 
$\iota'(c_{\frac12 - \delta,1}) \geq d-1$, since $L(c_{\frac12 - \delta,1}) \geq
\frac12$.  Therefore, we have
$$
 \iota(c) \geq \iota'(c) \geq
  \iota'(c_{0,\frac12 - \delta}) + \iota'(c_{\frac12 - \delta,1}) \geq
  1+d-1=d \,\,,
$$
which is a contradiction. $\Box$
\vspace{5mm}

As a consequence, we have

\vnn
{\bf Lemma 4.}\quad{\it
 Jacobi fields in the space ${\cal V}'_{c_{0,1}}$ are constant 
 multiples of the fields $\sin(ks)V(s)$; $0 \leq s \leq 1$ while 
 the negative eigenfields of $H_E$ are constant multiples of 
 the fields $\sin(\frac{k}{2}s) V(s)$, where 
 $V$ is any parallel section in the bundle $U(\perp c_\n{[0,1)})$ 
 (${\cal V}'_{c_{0,1}}$ being the space of all square integrable 
 normal vector fields along the arc $c_{0,1} := c\n_{[0,1]}$ 
 with Dirichlet boundary condition, as is defined in \S2).
}

\vnn
{\bf Remark.}\quad{\it A priori, the holonomy along the loop $c$ might be
non-trivial.  So, a parallel vector field $V$ of elements in $U(\perp
c\n_{[0,1)})$ might not close up at $s=1$.  Later we will show that the
holonomy along $c$ is trivial.}

\vnn
{\it Construction of the Araki Simplex.} \quad
 Let $V_1, \ldots ,V_{d-1}$ be parallel vector fields of orthonormal 
 elements in $U\perp c\n_{[0,1)}$.  By compactness argument, there 
 exists an $\eta > 0$ so that each orthogonally trajecting geodesics 
 $t \longmapsto \exp tx$ where $x \in U\perp c$ has no point focal 
 to $c$ in $t < \arctan \eta$.  Define $2(d-1)$ vector fields $X_i(s)$ 
 and $Y_i(s)$ ($0 \leq s \leq 1$) along $c$ as follows.  These vector 
 fields are not continuous at $s = 0$ and $s = \frac12$.
 \[
  X_i(s)=\left\{ \begin{array}{ll} V_i(s) & {\rm if\,\,}0 \leq s
   \leq \frac12\\ 0 & {\rm if\,\,}\frac12<s<1 \end{array}
   \right. 
 \]
 and
 \[
  Y_i(s)=\left\{ \begin{array}{ll}
   0 & {\rm if\,\,}0\leq s \leq \frac12\\
   V_i(s) & {\rm if\,\,}\frac12<s<1\,\,.
 \end{array}
 \right.
 \]
 Let $x = (x_1, \ldots, x_{d-1}) \in I \subset \eu^{d-1}$ and 
 $y = (y_1,\dots,y_{d-1}) \in I \subset \eu^{d-1}$, where $I$ 
 is a small interval in $\eu^{d-1}$ centered at the origin.  We define
 a $2(d-1)$-dimensional differentiable simplex $\wt{\sigma}$ 
 in $\Omega$ (here we regard $\Omega$ as a Riemannian Hilbert 
 manifold consisting of absolutely continuous alosed curves with 
 $L_2$ inner product) as follows
 $$
  \wt{\sigma}(x,y)(s) = \exp_{c(s)} \arctan \{ \eta\sin 
   (2 \pi s)(\sum_{i=1}^{d-1}(x_i X_i(s) + y_i Y_i(s))) \} \,\,.
 $$
 \noindent
 Here, by  ``arctan'' of a vector, we  will mean for a vector 
 $x \in U \perp c$ the resized vector $
  (\arctan \| x\|)\frac{x}{\|x\|}$.  W. Ballman pointed out to 
 us that Araki [A] constructed a simplex in the same way, i.e., 
 varying Jacobi fields independently outside the zero set, when 
 $M$ is a symmetric space.  So, such
 a simplex may be called {\it Araki simplex}.  
 The Araki simplex $\widetilde{\sigma}$ consists of curves all 
 passing through $c(0)$ and $c(\frac12)$.

\vnn
{\it Deformation of the Araki Simplex.}\quad
 We deform the Araki Simplex just constructed 
 in the following way.  If $x = y$ 
 we make  no change on the corresponding loop.  
If $x \not= y$, then we make suitable short cuts at the 
non-trivial angle created by the discrepancy $x \not= y$ 
at $s = \frac12$. 
For instance, we fix a small positive number $\delta$ and make a 
short cut between points corresponding to 
$s=\frac12 - \delta$ and $s = \frac12 + \delta$ 
by a small geodesic arc. 
After performing this modification and reparametrizing the 
corresponding loops by arc length, we get a $2(d-1)$-dimensional 
simplex $\sigma$ (we call this the ``short cut modification''). 
We note that
 \svnn
 
 (i) the intersection $\wt\sigma \cap \sigma$ consists 
  of those closed curves that are generated by $x=y$ where $x(=y)$ 
  satisfies the 
  condition that the parallel vector field $\sum_{i=1}^{d-1}x_iV_i$ 
  along $c$ closes 
  up at $s=1$, i.e., variations which integrate global Jacobi fields 
  on $c\n_{[0,1)}$, 
 \svnn

 and
 \svnn

 (ii) the vector 
  fields $\sin (2\pi s)X_i(s)$ and $\sin(2\pi s)Y_i(s)$ are naturally 
  regarded 
  as Jacobi fields along $c\n_{[0,\frac12]}$; respectively, 
  $c\n_{[\frac12,1]}$ which vanish at end points. 
  \svnn

 If $x \not= y$, then, after performing the 
 above modification, we see that $\sigma(x,y)$ is strictly under the level set 
 $\Omega^{=1}$ of  $E = 1$.  We see this, by applying Rauch type comparison 
 theorem of Berger (Rauch's second comparison; see [CE (\S 1.10)]) to 
 variations 
 $$
  \exp_{c(s)} \arctan \{ \sin (2\pi s) \sum_{i=1}^{d-1}
   x_iX_i(s)\} \,\,, \qquad s \in [0,\frac12] \,\,, \,\,
   t \in [0,\eta]
 $$ 
 of $c\n_{[0,\frac12]}$, and 
 $$
  \exp_{c(s)} \arctan \{ \sin (2\pi s) \sum_{i=1}^{d-1}
   x_iY_i(s)\} \,\,, \qquad s \in [\frac12,1] \,\,, \,\,
   t \in [0,\eta]
 $$ 
 of $c\n_{[\frac12,1]}$ in $\widetilde\sigma$.  Note that the 
 corresponding comparison variations in $S_{k}^{2}$ generate great 
 semicircles fixed at the north and south poles. 
\vnn

Summing up, we have
 
\vnn  
{\bf Lemma 5.}\quad{\it
 There exists a neighborhood $W$ of $c \in \cal C$ in $\Omega$ so 
 that the $2(d - 1)$-dimensional simplex $\sigma \cap W$ is 
 contained in $\Omega^1$.
}
\vspace{5mm}

We now exhibit the unstable simplex $\tau$ mentioned just after 
the statement of Lemma 2.  In fact, $\tau$ is 
the $(d - 1)-$dimensional subsimplex of the modified Araki simplex 
$\sigma$ corresponding to the parameters $x = - y$.  The simplex 
$\tau$ is not itself the strong unstable simplex.  However, because $\tau$ 
and the strong unstable simplex $\tau'$ constructed by exponentiating 
the negative eigenspace of $H_E$ are, downstairs in $M$, 
both contained in a tubular neighborhood of the geodesic $c$, we see 

\vnn
{\bf Lemma 6.}\quad{\it There exists a neighborhood $U$ of $c$ in $\Omega$ 
such that, for $\ep>0$ sufficiently small and a subneighborhood $W\subset U$, 
$\tau$ constructed above represents the same homotopy class as $\tau'$ in 
$\pi_{d-1}(W,W\cap \Omega^{1-\ep})$.
}
\vspace{5mm}

{\it Proof.} Both strong unstable simplex $\tau'$ and the simplex $\tau$ 
are constructed by exponentiating certain variation vector field in a negative 
cone in $T_{\wt c}\Omega$ of the $H_E$. 
Moreover, choosing a neighborhood $U$ of $c$ in $\Omega$ sufficiently small, 
we may assume that the set of critical points of $E$ in $U$ coincides 
with the connected critical submanifold, say, $C$, containing $c$ which, 
at each point $\wt c\in C$,  is tangent to $\cal Z$, where $\cal Z$ is 
the zero eigenspace in $T_{\wt c}\Omega$ of $H_E$ (ref. Theorem E). 
Therefore, if we choose sufficiently small subneighborhood $W\subset U$ 
of $c$ (and therefore $\ep>0$ sufficiently small), we see that there exists a 
homotopy connecting $\tau$ and $\tau'$ in the space 
$\pi_{d-1}(W,W\cap\Omega^{1-\ep})$. 
$\Box$
\vnn

Lemma 6 implies :

\vnn 
{\bf Lemma 7.}\quad{\it If $c \in {\cal C}^*$, then there is a neighborhood 
$U$ of $c$ in $\Omega$ so that, for $\ep>0$ sufficiently small 
and a subneighborhood $W \subset U$, $\tau$ constructed inside the modified 
Araki simplex represents a nontrivial element 
in $\pi_{d-1}(W,W \cap \Omega^{1-\ep})$.}
\vnn
\indent
The reason why we have chosen the ``short cut  modification" is explained 
in the following way. 
If we define a $(d-1)$-dimensional simplex
$\overline \tau$ by 
$$\overline\tau(x):=\exp_{c(s)}\arctan \{\eta\sin(\pi
s)\sum_{i=1}^{d-1}x_iV_i(s)\}\,\,,$$
\noindent then direct calculation of the
Hessian implies that $\overline\tau$ also defines an unstable simplex at $c$
and belongs to the same class as $\tau$ in the relative homotopy group.  This
unstable simplex corresponds to the eigen vector of the index form with
negative eigenvalue.  In this sense, $\overline \tau$ is more natural than
$\tau$.  Now define a $(2d-2)$-simplex $\overline\sigma$, which also contains
the $(d-1)$-simplex (the one defined by $x=y$ in our simplex $\sigma$) 
corresponding to the global Jacobi field on $c$, by 
$$\overline\sigma(x,y)(s):=\exp_{c(s)}\arctan\{\eta\sum_{i=1}^{d-1}(x_i\sin(\pi
s)+y_i\sin(ks))V_i(s)\}\quad (k=2\pi)\,\,.$$
\noindent
Although this construction is natural, it turns out that it is not clear
whether there exists an interval $I$ containing $0$ such that $\overline\sigma(I
\times I)$ is contained in $\Omega^1$.  This is the reason why our
construction of $(2d-2)$-simplex $\sigma$ is based on the short cut argument
of broken geodesics in the model space, although the unstable simplex $\tau$
does not directly integrate the negative eigenspace of the Hessian of the
energy functional $E$ at $c$.\vnn
\indent

We return to our ``short cut modification" and consider the holonomy problem 
mentioned just after Lemma 4. 
One of the following two cases is possible.  Namely, either
\vnn
(A) For at least one choice of $x_0 \in I$, there is some $\ep>0$ such that
$$\begin{array}{c}
\exp_{c(s)}\arctan\{t\sin(2\pi
s)(\sum_{i=1}^{d-1}x_{0,i}X_i(s)+x_{0,i}Y_i(s))\}\\[5mm]
=\exp_{c(s)}\arctan\{t\sin(2\pi s)(\sum_{i=1}^{d-1}x_{0,i}V_i(s))\}\,\,,
\end{array}$$
\noindent
is contained in $\Omega^{1-2\ep}$ for
all $t \in (\eta/2,\eta]$. In the picture of this situation, 
we find two variation vector fields $V$ and $Y$ 
along $c$ of the form  
$$\begin{array}{ccl} Y & = & \sin(2\pi
s)\sum_{i=1}^{d-1}x_iV_i(s)\,\,\,\,\,\,({\rm Jacobi\,\,fields})\,\,,\\[5mm]
V & = & \sin(2\pi s)
\sum_{i=1}^{d-1}(x_iX_i(s)-x_iY_i(s))\,\,\,({\rm tangent\,to}\,\tau)\,\,,
\end{array}$$
\noindent
outside a small neighborhood of $s=\frac12$. 
In the picture of the simplex $\sigma$, we find a $(d-1)$-dimensional 
simplex $\tau\cap W$ which lines in $\Omega^{<1}$ except at $c$ 
and moreover we have extra one direction represented by $x_0$ which 
also behaves exactly like a strong unstable simplex, 
\vnn

or else,
\vnn
(B) There exists an $\alpha$; $0<\alpha<1$ so that whenever 
$\n x_1\n,\dots,\n x_{d-1}\n \leq \alpha$,
$$\exp_{c(s)}\arctan\{\sin(2\pi s)(\sum_{i=1}^{d-1}x_iV_i(s))\}$$
\noindent
lies in $\Omega^{=1}$.  (In particular, each
parallel vector field 
$\sum_{i=1}^{d-1}x_iV_i(s)$ closes up at $s=1$.)
\vnn
\indent 
If Case (A) prevails, $\tau \cap W$ rides on a $d$-dimensional
submanifold of $\sigma \cap W$ which lies in $\Omega^{<1}$ except at $c$ and
hence $\tau \cap W$ can be deformed into $W \cap \Omega^{1-\ep}$.  This
contradicts the conclusion of Lemma 7.  Hence $c \not\in {\cal C}^*$.  If, on
the other hand, we start out with a $c \in {\cal C}^*$, then Case (B) must
really be the case.  In particular, the holonomy along $c \in {\cal C}^*$
must be trivial.  We thus get a $(d-1)$-dimensional local submanifold $S$ of 
$\Omega^{=1}$ which is tangent to the 0 eigenspace of the Hessian of $E$ 
defined on ${\cal V}_{c_{0,1}}$ through $c$.
\vnn 
{\bf Lemma 8.}\quad{\it In the present
situation, each parallel vector field
$\sum_{i=1}^{d-1}x_iV_i(s)$ closes up at $s=1$, i.e., the holonomy along $c$
is trivial, and each member $\wt c$ of $S$ is a (smooth) closed geodesic in
$\cal C$.}
\vnn
\indent
{\it Proof.} We need to prove the second assertion.  If $\wt
c$ is not a critical point of $E$, there exists at least one $x_0 \in I
\subset \eu^{d-1}$ such that the $(d-1)$-dimensional simplex defined by the
$(d-1)$-dimensional affine subspace through $x_0$ orthogonal to the linear
subspace defined by $x=y$ contains no critical point.  Then, by following the
trajectory of $-{\rm grad}\,\,E$, $\tau$ is deformed into $W \cap
\Omega^{1-\ep}$, which contradicts the assumption that we started with $c \in
{\cal C}^*$.  Hence all $\wt c \in S$ are closed geodesics.  If some $\wt c
\not\in \cal C$, then, it follows from Lemma 3 and its proof that $\iota(\wt c)>d-1$. 
So $\tau$ is again deformed into $W
\cap \Omega^{1-\ep}$ via the unstable simplex of $\wt c$.  Hence, either way,
we get a contradiction.  $\Box$\vnn \indent By construction, we also see that
for any $\wt c \in S$, $\wt c(0)=c(0)$.  Translated into $M$, this means that
there is an open tube $B$ (cone-like at $s=0$ and $\frac12$) 
around the set $c(0,\frac12) \cup c(\frac12,1)$
such that for each $q \in B$, a geodesic joining $c(0)$ to $q$ extends to a
closed geodesic in $\cal C$ whose image lies in $B$ except at $s\in \frac12\z$.  
Applying Lemma 3 to each geodesics proves Lemma 2. 
$\Box$
\vnn
\indent
Even more is true.  
\vnn {\bf
Lemma 9.}\quad{\it Let $c \in {\cal C}^*$ and let $\cal U$$\subset T_{c(0)}M$
be the set in Lemma 2.  Then, there exists an open set ${\cal U}^*$; $c'(0)
\in {\cal U}^* \subset \cal U$, so that, for all $u \in {\cal U}^*$, $c_u \in
{\cal C}^*$.}\vnn
\indent
{\it Proof}. Since $c\in \cal C^*$ is a closed geodesic in the compact Riemannian 
manifold $M$, there exists an $\ep>0$ such that there are no critical values 
for $E$ on $\Omega$ in the intervals $(1-\ep,1)$ and $(1,1+\ep)$. Lemma 8 
implies that there exists an open set $\cal U^*$ in $T_{c(0)}M$ containing 
$c'(0)$ satisfying the condition that $u\in\cal U*$ implies $c_u\in \cal C$. 
The strong unstable simplex constructed by exponentiating the negative 
eigenspace of $H_E$ defines a differentiable simplex 
$\tau:(I,\partial I) \longrightarrow (\Omega^{1+\ep},\Omega^{1-\ep})$ 
($I\subset \eu^{d-1}$). We introduce the he compact-open topology to the 
set $\Sigma$ of all absolutely continuous maps 
$(I,\partial I)\longrightarrow  (\Omega^{1+\ep},\omega^{1-\ep})$, by which 
we can argue the closeness of maps in $\Sigma$. 
Then if $\cal U^*$ is a sufficiently small open set containing $c'(0)$ 
in $T_{c(0)}M$, then, for $\forall u\in\cal U^*$, 
we can construct, by exponentiating $(d-1)$-dimensional negative eigenspace 
of $H_E$ at $c_u$, a $(d-1)$-dimensional strong unstable simplex 
which is homotopic in $\Sigma$ to $\tau$ constructed above. $\Box$
\vnn

That is to say, the set
\begin{center}
${\cal U}^*=\{u \in UT_{c(0)}M\,;\,c_u \in {\cal C}^*\}$
\end{center}
is an open set in $UT_{c(0)}M$.  On the other hand, by Lemma 1 and the
continuous dependence of geodesics on their initial values, the set
${\cal U}^*$ is also a closed set.  Since $UT_{c(0)}M$ is connected, ${\cal U}^*$
must in fact be all of $UT_{c(0)}M$.  Together with Lemma 2, we summarize our
result as
\vnn
{\bf Lemma 10.}\quad{\it Let $M$ be assumed in this section.  Then, there exists
a point $p \in M$ such that for all $u \in T_pM$, $c_u$ is a closed geodesic of
prime length $1$ and $K(\tau)=k^2$ for all 2-planes $\tau$ tangent to the radial
direction from $p$.}\vnn
\indent
{\it Proof.} Take a $c \in {\cal C}^*$ and let $p:=c(0)$.  $\Box$\vnn
\indent
Now it is a standard technique to construct an explicit isometry from $M$ onto
$S^d_k$ exactly as in Toponogov's Maximum diameter theorem (see, for instance, 
[CE (\S 6.5)] or [GKM (\S 7.3)]).  Thus, the Main Theorem is proven as soon as
Lemma 1 is established.

\section{Proof of Lemma 1}

In this section we work on the finite-dimensional approximation (because 
we need some analytic argument which is described simpler 
in the finite-dimensional approximation). 
Recall that $\Omega$ is the space of all absolutely continuous closed curves 
with $L^2$ derivative and in particular contains piecewise differentiable curves. 
It is well-known that for each $r > 0$, $\Omega^r$ 
(it is defined just after Theorem D in \S2) contains a
submanifold $'\Omega_r$ which is diffeomorphic to an open set 
in some finite product $M \times \cdots \times M$ and homotopy 
equivalent to $\Omega^{<r}$.  $'\Omega='\Omega_r$ consists of 
broken geodesics. 
The functional $E$ becomes a proper 
function on $'\Omega_r$.  The space $'\Omega_r$ contains all the 
critical points in $\Omega^{<r}$ and the Hessian of 
$E\n_{'\Omega_r}$ retains the same index as $E$ at each critical 
point.  
Moreover, Theorem E in \S2 remains true in this finite-dimensional setting. 
For details, see [Mi (\S 16)] and [Bo].  If $a<r$, we put 
$'\Omega^a_r := \, '\Omega_r \cap \Omega^a$.

In this section, we continue to assume $K \geq 4\pi^2$.  The following proposition 
is essentially contained in some earlier works of M. Berger and is easy to 
prove by Morse-Schoenberg index comparison with $S_k^d$ and the tautological 
isomorphism $\pi_i(\Omega) \cong \pi_{i+1}(M)$.
\vnn
{\bf Proposition1.}\quad{\it If $M$ contains no closed geodesic of length $\leq 1/2$, 
then $M$has the homotopy type of a sphere.  
In particular, we have} 
\[ \pi_i(\Omega,M)
\cong \left \{ 
\begin{array}{ll} \z & {\rm if\,\,}i=d-1\\
0 & {\rm for\,\,}0 \leq i \leq d-2
\end{array}
\right. \]
\noindent
for the relative homotopy groups $\pi_i(\Omega,M)$ up to $i \leq d-1$.
\vnn
\indent
We now return to the assumption that the length of the shortest closed
geodesics on $M$ is 1.  Let $\cal C$ and ${\cal C}^*$ be as defined in \S 3. 
We wish to prove that a strong unstable simplex at at least one $c \in \cal
C$ represents a nontrivial class of $\pi_{d-1}(\Omega,M)$.  Our technique
will be to approximate $E$ with other functionals that are guaranteed to have
nontrivial unstable simplexes.  Although all our arguments carry through in
all of $\Omega$ in an $S^1$-invariant fashion, essentially because the
functional $E$ satisfies the Condition (C) of Palais and Smale and because an
$S^1$-invariant formulation of Theorem E is available [GM2], we find it a
little easier to work in a finite dimensional space.\nln
\indent
More precisely, choose $r$ sufficiently large, say $r>2$.  Then, all closed
geodesics not in $\Omega^r$ will have index $>2(d-1)$.  Let
$'\Omega:=\,'\Omega_r$.  Then $$\pi_i('\Omega,M) \cong \pi_i(\Omega,M)$$
\noindent
for all $i$; $0 \leq i \leq 2d-3$, and $d-1 \leq 2d-3$ if $d \geq 2$.  Using
Theorem E and a partition of unity on $'\Omega$, we can approximate $E$ with a
sequence $\{E_n\}_{n=1}^{\infty}$ of functionals on $'\Omega$ with the
following properties.\svnn

(i) $\lim_{n \rightarrow \infty} E_n=E$ in the $C^2$
topology.\svnn

(ii) For some $\ep>0$, all critical points of $E_n$ in the
closure of the set $L:=\,{'\Omega_{1+\ep}}-\,{'\Omega_{1-\ep}}$ either belong
to $\Omega_{=1}$ or have index $\geq 2d-2$, and outside $L$, each $E_n$ agrees 
with $E$.
\svnn

(iii) Each $E_n$ has only nondegenerate critical points in the
set $L$, all of which have index $\geq d-1$.
\svnn

Let $C$ be the set of all closed geodesics in $'\Omega^{=1}$ and let $C_n$ be
the set of all critical points of $E_n$ that lie in $L$.
\vnn
{\bf Lemma 11.}\quad{\it For each $n$, there exists in $C_n$, at least one
critical point of $E_n$ that possesses a strong unstable simplex that represents
a nontrivial element in $\pi_{d-1}('\Omega,M)$.}
\vnn
\indent
{\it Proof.} From the topology described in Proposition 1, there must exist a
nontrivial element $\rho$ of $\pi_{d-1}('\Omega,M)$.  We first deform $\rho$ so
that the only points of $C_n - (M \cap C_n)$ that lies on the image of $\rho$
are the relative maxima of $E_n \circ \rho$.  In fact, since there are no
critical points of index $<d-1$ except in $M$, at every critical point of $E_n$
lying on $\rho$, say $c$, other than relative maxima, the unstable dimension of
$E_n$ in $'\Omega$ is strictly greater than the unstable dimension of
$E_n \circ \rho$ in the image of $\rho$.  Therefore, in some neighborhood of $c$
in which a chart of the form described in Theorem E is valid, we can deform
$\rho$ in a direction transversal to itself and which decreases $E_n$.  Since
the critical points of $E_n$ are isolated and $\rho$ is contained in a compact
region, by repeating this deformation a finite number of times and by deforming
$\rho$ along the trajectory of $-{\rm grad}\,\,E_n$, we can deform $\rho$ until
it is expressed as a sum of disjoint simplexes, each summand of which is a
simplex in $('\Omega,M)$, hanging from a single critical point of index $=d-1$.
Such critical points must be in $C_n$, and at least one summand must be
nontrivial itself. We deform this simplex by a differentiable homotopy, 
if necessary, into a strong unstable simplex (in the sense defined just 
after Theorem E) in $\pi_{d-1}('\Omega,M)$. $\Box$
\vnn
\indent
Of course, it is not necessarily true that a sequence of critical points
$\{c_n\}$ of $C_n$ converges to a closed geodesic.  However, that $\lim_{n
\rightarrow \infty} C_n \subset C$ in the following weaker sense is clear.
\vnn
{\bf Lemma 12.}\quad{\it Given any open neighborhood $\cal U$ of $C$ in
$'\Omega$, whenever $n$ is large enough, $C_n \subset \cal U$.}
\vnn
\indent
In fact, since the convergence is specified in the $C^2$ topology, we can state
the even stronger
\vnn
{\bf Lemma 13.}\quad{\it Let $\{ {\cal U}^-_c \subset {\cal U}^{-0}_c\}_{c
\in C}$ be a family of pairs of open sets in $'\Omega$ so that, for each
$c \in C$, ${\cal U}^-_c$ is a neighborhood of the strong unstable
submanifold $U^{-}_c$ of $E$ at $c$ and ${\cal U}^{-0}_c$ is a neighborhood of
the unstable submanifold $U^{-0}_c$.  Then, for $n$ sufficiently large, for
each $c_n \in C_n$, there exists some $c \in C$, so that $U_{c_n}^{-}$, the
strong unstable manifold of $E_n$ at $c_n$ is contained in ${\cal U}^{-0}_c$. 
Moreover, for such $c_n$ and $c$, a strong unstable simplex $\tau_n$ of $c_n$
contains a subsimplex $\tau'_n$ with $\dim \tau'_n=\dim U_c^-=\iota(c)$ which
is actually contained in ${\cal U}^-_c$.}
\vnn
\indent
To see the above, we can take a local coordinate expression around each $c \in
C$ as described in Theorem E and look at the partial derivatives.  By taking
$n$ large, if $c_n \in C_n$ is close to $c \in C$, the corresponding second
derivatives respectively of $E_n$ at $c_n$, $E$ at $c_n$ and $E$ at $c$ can all
be made arbitratily close to each other by the property (i).  But, in $\cal U$, the
strong unstable submanifolds and unstable submanifolds are determined by the
second partial derivatives.
\nln
\indent
Now, for each $n$, let $c_n$ be the critical point in Lemma 11 which has a
strong unstable simplex $\tau_n$ that is nontrivial in $\pi_{d-1}(\Omega,M)$. 
For such a $c_n$, $\tau_n \cap \cal U$ must itself be contained in a neighborhood 
${\cal U}^-_c$ of the strong unstable submanifold at some $c \in C$ by index 
comparison and the dimensional consideration.  From the construction of 
$\tau_n$, this $c$ must be $\in \cal C$.  Let $\tau$ be a strong unstable 
simplex at $c$ with $\tau(\partial I^{d-1}) \subset M$.  By repeating the 
standard Morse theoretic arguments as in Lemmas 6-7, we see that $\tau$ 
represents a nontrivial element in $\pi_{d-1}(\Omega,M)$.  Hence $c \in 
{\cal C}^*$.  Then, that ${\cal C}^*$ is closed follows from the topological 
arguments in the proof of Lemma 9. 
This completes the proof of Lemma 1 and thus of Main Theorem.  
$\Box$

\vspace{5mm}
\noindent
Y. I.
\nln
Department of Information and Communication Engineering\nln
Fukuoka Institute of Technology
\nln
R. K.
\nln
Graduate School of Mathematics
\nln
Nagoya University
\end{document}